\documentclass[10pt, reqno]{amsart}
\usepackage{amssymb,latexsym}
\usepackage{eucal}

\newtheorem*{proposition*}{}

\newtheorem{lemma}{Lemma}

\newtheorem*{Lem17.3}{Lemma 17.3}
\newtheorem*{Lem18.3}{Lemma 18.3}
\newtheorem*{Lem18.5}{Lemma 18.5}
\newtheorem*{Lem19.0}{Lemma 19.0}
\newtheorem*{Lem19.1}{Lemma 19.1}
\newtheorem*{Lem19.1.3}{Lemma 19.1.3}

\newtheorem*{Th}{Theorem}
\newtheorem*{Cor}{Corollary}
\newtheorem*{Pr1}{Problem 1}
\newtheorem*{Pr2}{Problem 2}
\newtheorem*{Pr3}{Problem 3}

\newcommand{\la}{{\langle}}
\newcommand{\ra}{{\rangle}}
\newcommand{\e}{{\varepsilon}}

\newcommand{\R}{\mathcal{R}}  \newcommand{\Ss}{\mathcal{S}}

\newcommand{\D}{\Delta}
\newcommand{\A}{\mathcal{A}}

\newcommand{\X}{\mathcal{X}}

\newcommand{\p}{\partial}

\newcommand{\C}{\mathcal{C}}

\newcommand{\MM}{\mathfrak{M}}
\newcommand{\wtl}{\widetilde}

\begin{document}

\title[On varieties of groups in which all periodic groups are abelian]
{On varieties of groups in which all periodic groups are abelian}
\author{S.V. Ivanov}
\address{Department of Mathematics\\
University of Illinois \\
Urbana,  IL 61801, USA}
\email{ivanov@math.uiuc.edu}
\thanks{The first author is supported in part by NSF grant DMS 00-99612}

\author{A.M. Storozhev}
\address{Australian Mathematics Trust\\
University of Canberra \\
Belconnen, ACT 2616, Australia}
\email{andreis@amt.canberra.edu.au }


\begin{abstract}
To solve a number of problems on varieties of groups, stated by
Kleiman, Kuznetsov, Ol'shanskii, Shmel'kin in the 1970's and
1980's, we construct continuously many varieties of groups in
which all periodic groups are abelian and  whose pairwise
intersections are the variety of all abelian groups.
\end{abstract}
\maketitle

\section{Introduction}

The first example of a nonabelian variety of groups in which all
finite groups are abelian was constructed by Ol'shanskii
\cite{O85}. This example provided a positive solution to a problem
of Hanna Neumann \cite[Problem 5]{N67}  on the existence of such
varieties of group and demonstrated that a nonabelian variety of
groups might have a very limited intersection with the class of
finite groups. Shmel'kin \cite[Problem 4.73(b)]{KN82} posed a
strengthened version of this Hanna Neumann's problem by  asking
about the existence of a nonabelian variety of groups in which all
periodic groups are abelian. A positive solution to this problem
of Shmel'kin was announced by the authors in \cite{IO91},
\cite{S91}. The aim of this article is to present details of the
following construction which is sketched in \cite[Theorem 3]{IO91}
and which solves Shmel'kin's problem (together with a number of
other problems mentioned below; for the sake of simplicity of
proofs we change the identities of \cite{IO91}).

Let $\Omega =(\Omega(1), \Omega(2),\dots)$ be an infinite
sequence of numbers $\Omega(k)\in \{ 0, 1\}$, $k = 1,2, \dots$,
let  $p_k$ be the $k$th prime number, $q_k =(p_1 \dots p_k)^k$,
and set
\begin{equation}
v_k(x,y) =[[x^d, y^d]^d, x^{d q_k }] , \label{D:1}
\end{equation}
\begin{multline}
w_{\Omega, k}(x,y) = [x,y]^{\Omega(k) + \e_1} v_k(x,y)^{d^k
n+1} [x,y]^{\e_2} v_k(x,y)^{d^k n+2} \dots  \\ \dots
[x,y]^{\e_{h-1}} v_k(x,y)^{d^k n+h-1} [x,y]^{\e_{h}}
v_k(x,y)^{d^k n+h}, \label{D:2}
\end{multline}
where $[a,b] = aba^{-1}b^{-1}$ is the commutator of $a$ and
$b$, $h \equiv 0 \pmod{10}$,
\begin{gather*}
\e_{10\ell+1}=  \e_{10\ell+2}=  \e_{10\ell+3}=  \e_{10\ell+5}=
  \e_{10\ell+6}=  1, \\
\e_{10\ell+4}=  \e_{10\ell+7}=  \e_{10\ell+8}=  \e_{10\ell+9}=
\e_{10\ell+10}= -1 ,
\end{gather*}
$\ell=  0,1,\dots ,h/10-1$, and $h, d, n$ are sufficiently
large positive integers (with $n \gg d \gg h \gg 1$). The
following is our  main result.

\begin{Th}[]
\label{Th} Let $\Omega = ( \Omega(1), \Omega(2), \dots )$ be an
infinite sequence of numbers $\Omega(1)$, $\Omega(2), \ldots \in
\{0, 1\}$ that contains infinitely many $1$'s and $\MM_{\Omega}$
be the variety of groups defined by identities $w_{\Omega, k}(x,y)
\equiv 1$, $k= 1,2,\dots$, where the words $w_{\Omega, k}(x,y)$
are given by formula \eqref{D:2}. Then $\MM_{\Omega }$ is a
nonabelian variety of groups in which all periodic groups are
abelian. Furthermore, if $\Omega_1 \neq \Omega_2$ then the
intersection $\MM_{\Omega_1} \cap \MM_{\Omega_2}$ is the variety
$\mathfrak A$ of all abelian groups.
\end{Th}

Let $\mathfrak {U}$ be a variety of groups. Recall that a variety
of groups $\mathfrak V$ is called   {\em just-non-$\mathfrak U$}
if $\mathfrak V$ properly contains the variety $\mathfrak U$, i.e.
$\mathfrak U \subset \mathfrak V$, and there is no group variety
$\mathfrak V'$ such that $\mathfrak U \subset \mathfrak V' \subset
\mathfrak V$. It follows from Zorn's lemma that every group
variety $\MM_{\Omega}$ of our Theorem contains a
just-non-$\mathfrak A$ variety $\widetilde \MM_{\Omega}$, where
$\mathfrak A$ is the variety of all abelian groups. Since the set
of sequences $\Omega$ with infinitely many $1$'s is continuous and
$\widetilde \MM_{\Omega_1} \cap \widetilde \MM_{\Omega_2} =
\mathfrak A$ if $\Omega_1 \neq \Omega_2$, we have the following.

\begin{Cor}[] Let $\mathfrak A$ denote the variety of all abelian groups.
Then the set of just-non-$\mathfrak A$ varieties of groups is
continuous.
\end{Cor}

We remark that Kleiman \cite{K83} earlier found a solvable variety
of groups $\mathfrak K$ such that the set of  just-non-$\mathfrak
K$ varieties is continuous. This Corollary provides a new insight
into the structure of the lattice of group varieties and enables
us to solve the following three problems posed in the 1970's and
1980's.

\begin{Pr1}[Ol'shanskii, Problem 4.46(b) \cite{KN82}]
How many varieties are there that have no finite basis for
their identities and in which any proper subvariety has a
finite basis for its identities?
\end{Pr1}

\begin{Pr2}[Kuznetsov, Problem 6.16 \cite{KN82}]
How many varieties are there that contain only finitely many or
countably many  subvarieties?
\end{Pr2}

\begin{Pr3}[Kleiman, Problem 8.20 \cite{KN82}]
How many  just-non-abelian (or  just-non-nilpotent) varieties
of groups are there?
\end{Pr3}

The authors are grateful to A.Yu.~Ol'shanskii for pointing out
that Kozhevnikov \cite{Ko00a}, \cite{Ko00b} earlier constructed
continuously many just-non-$\mathfrak A_n$ varieties of groups,
where $n \gg 1$ is odd and $\mathfrak A_n$ is the variety of all
abelian groups of exponent $n$. In particular, this result of
Kozhevnikov also implies that the sets of varieties of groups in
Problems 1--3 are continuous. The authors also wish to thank the
referee for useful remarks.

\section{Proof of Theorem}

It is fairly  easy to see that all periodic groups in the variety
$\MM_\Omega$ of our Theorem  are abelian. Indeed, it follows from
the identity $w_{\Omega, k}(x,y) \equiv 1$, where $w_{\Omega,
k}(x,y)$ is given by \eqref{D:2}, that if $\Omega(k) =1$ then the
quasiidentity $x^{q_k} =1 \to [x,y] =1$ holds in the variety
$\MM_\Omega$. In particular, any element of finite order of a
group $G \in \MM_\Omega$ lies in the center of $G$. It is also
clear that if $\Omega_1 \neq \Omega_2$, say $\Omega_1(k_0) \neq
\Omega_2(k_0)$ for some $k_0 $, then it follows from identities
$w_{\Omega_1, k_0}(x,y) \equiv 1$, $w_{\Omega_2, k_0}(x,y) \equiv
1$ that $[x,y]\equiv 1$ which means that $\MM_{\Omega_1} \cap
\MM_{\Omega_2} = \mathfrak A$. Therefore, the only nontrivial part
of our Theorem is to show that $\MM_{\Omega}$ is a nonabelian
variety of groups.

\medskip

To prove that $\MM_{\Omega}$ is  nonabelian, we will construct
a presentation for a free group of rank $m > 1$ in $\MM_\Omega$
by means of generators and defining relations and use the
geometric machinery of graded diagrams,  developed by
Ol'shanskii \cite{O85}, \cite{O89}, to study this group. In
particular, we will use the notation and terminology of
\cite{O89} and all notions that are not defined in this paper
can be found in \cite{O89}.

As in \cite{O89},  we will use numerical parameters
$$
\alpha \succ \beta  \succ \gamma  \succ \delta \succ \e \succ
\zeta \succ \eta \succ \iota
$$
and $h = \delta^{-1}$, $d = \eta^{-1}$, $n = \iota^{-1}$ ($h,
d, n$ were already used in \eqref{D:1}--\eqref{D:2}) and employ
the least parameter principle (LPP) (according to LPP a small
positive value for, say,  $\zeta$ is chosen to satisfy all
inequalities whose smallest (in terms of the relation $\succ$)
parameter is $\zeta$).

Let $\A = \{ a_1, \dots, a_m \}$  be an alphabet, $m > 1$, and
$F(\A)$ be the free group in $\A$. Elements of $F(\A)$ are
referred to as words in $\A^{\pm 1} = \A \cup \A^{-1}$ or just
words. Denote $G(0) = F(\A)$ and let the set $\R_0$ be empty.
Now consider an arbitrary infinite sequence $\Omega
=(\Omega(1), \Omega(2), \dots)$, where $\Omega(k) \in \{
0,1\}$, $k=1,2,\dots$. To define the group $G(i)$ by induction
on $i \ge 1$ for this given $\Omega$, assume that the group
$G(i-1)$ is already constructed by its presentation
$$
G(i-1) = \la \A \; \| \; R=1, R \in \R_{i-1} \ra .
$$
Let $\X_i$ be a  set of words (in $\A^{\pm 1}$) of length $i$,
called {\em periods of rank $i$}, which is maximal with respect
to the following two properties:
\begin{enumerate}
\item[(A1)] If $A \in \X_i$ then $A$ (that is, the image of $A$
in $G(i-1)$) is not conjugate in $G(i-1)$ to a power of a word
of length $< |A| = i$.

\item[(A2)] If $A$, $B$ are distinct elements of $\X_i$ then
$A$ is not conjugate in $G(i-1)$ to $B$ or $B^{-1}$.
\end{enumerate}

If the images of two words $X$, $Y$ are equal in the group
$G(i-1)$, $i \ge 1$, then we will say that $X$  is {\em equal
in rank} $i-1$ to $Y$ and write $ X \overset {i-1} = Y$.
Analogously, we will say that two words $X$, $Y$ are {\em
conjugate in rank} $i-1$ if their images are conjugate in the
group $G(i-1)$. As in  \cite{O89}, a word $A$ is called {\em
simple} in rank $i-1$, $i \ge 1$, if $A$ is conjugate in rank
$i-1$ neither to a power $B^\ell$, where $|B| <|A|$,  nor to a
power of period of some rank $\le i-1$. We will also say that
two pairs $(X_1, X_2)$, $(Y_1, Y_2)$ of words are conjugate in
rank $i-1$, $i \ge 1$, if there is a word $W$ such that  $X_1
\overset {i-1} = W Y_1 W^{-1}$ and $X_2 \overset {i-1} = W Y_2
W^{-1}$.

Consider the set of all possible pairs $(X, Y)$ of words in
$\A^{\pm 1}$ and pick a positive integer $k$. This set is
partitioned by equivalence $k$-classes $\C_\ell(k)$, $\ell=
1,2,\dots$,  of the equivalence relation $\overset { k} \sim$
defined by $(X_1, Y_1) \overset { k} \sim (X_2, Y_2)$ if and
only if the pairs $( v_{k}(X_1, Y_1),  w_{\Omega, k}(X_1, Y_1)
)$ and $( v_{k}(X_2, Y_2), w_{\Omega, k}(X_2, Y_2))$  are
conjugate in rank $i-1$. It is convenient to enumerate (in some
way)
$$
\C_{A^f, 1}(k) , \C_{A^f, 2}(k), \dots
$$
all $k$-classes of pairs $(X, Y)$ such that $w_{\Omega, k}(X,
Y) \overset {i-1} \neq 1$ and  $v_{k}(X, Y)$ is conjugate in
rank $i-1$ to some power $A^f$, where $A \in \X_i$ and $f =
f(X,Y)$ are fixed.

It follows from definitions that every class $\C_{A^f, j}(k)$
contains a pair
$$
(X_{A^f, j, k}, \bar Y_{A^f, j, k} )
$$
with the following properties. The word $X_{A^f, j, k}$ is
graphically (that is, letter-by-letter) equal to a power of
$B_{A^f, j, k}$, where $ B_{A^f, j, k}$ is simple in rank $i-1$
or a period of rank $\le i-1$; $\bar Y_{A^f, j, k} \equiv
Z_{A^f, j, k} Y_{A^f, j, k} Z_{A^f, j, k}^{-1}$, where the
symbol \ '$\equiv$' \  means the graphical equality, $Y_{A^f,
j, k}$ is graphically equal to a power of $C_{A^f, j, k}$,
where $C_{A^f, j, k}$ is simple in rank $i-1$ or a period of
rank $\le i-1$. We can also assume that if $D_1 \in \{ A,
B_{A^f, j, k}, C_{A^f, j, k} \}$  is conjugate in rank $i-1$ to
$D_2^{\pm 1}$, where $D_2 \in \{ A, B_{A^f, j, k}, C_{A^f, j,
k} \}$, then $D_1 \equiv D_2$. Finally, the word $Z_{A^f, j,
k}$ is picked for fixed $X_{A^f, j, k}$, $Y_{A^f, j, k}$ so
that the length $|Z_{A^f, j, k}|$ is minimal (and the pair
$(X_{A^f, j, k}, Z_{A^f, j, k} Y_{A^f, j, k} Z_{A^f, j,
k}^{-1})$ belongs to $\C_{A^f, j}(k)$). Similar to \cite{O85},
\cite{O89}, the triple $(X_{A^f, j, k}, Y_{A^f, j, k}, Z_{A^f,
j, k} )$ is called an $(A^f, j, k)$-{\em triple} corresponding
to the class $\C_{A^f, j}(k)$ (in rank $i-1$).

Now for every class $\C_{A^f, j}(k)$  we pick a corresponding
$(A^f, j, k)$-{\em triple}
$$
(X_{A^f, j, k},  Y_{A^f, j, k}, Z_{A^f, j, k} )
$$
in rank $i-1$ and define a defining word $R_{A^f, j, k}$ of
rank $i$ as follows. Pick a word $W_{A^f, j, k}$ of minimal
length so that
$$
v_{k}( X_{A^f, j, k},  \bar Y_{A^f, j, k} ) \overset {i-1} =
W_{A^f, j, k} A^{f} W_{A^f, j, k}^{-1} .
$$

Let   $T_{A^f, j, k}$, $U_{A^f, j, k}$   be  words of minimal
length such that
\begin{gather*}
T_{A^f, j, k}  \overset {i-1} = W_{A^f, j, k}^{-1}  [ X_{A^f,
j, k},  \bar Y_{A^f, j, k} ] W_{A^f, j, k}, \\
U_{A^f, j, k}  \overset {i-1} =
W_{A^f, j, k}^{-1}  [ X_{A^f, j, k},
\bar Y_{A^f, j, k} ]^2 W_{A^f, j, k}.
\end{gather*}

In accordance with (\ref{D:2}), if $\Omega(k) =0$, then we set
\begin{equation}
R_{A^f, j, k}  =  T_{A^f, j, k}^{\e_1} A^{(d^k n+1)f} T_{A^f,
j, k}^{\e_2} A^{(d^k n+2)f} \dots T_{A^f, j, k}^{\e_h} A^{(d^k
n+h)f}, \label{D:3}
\end{equation}
and, if $\Omega(k) =1$, then we set
\begin{equation}
R_{A^f, j, k}  =  U_{A^f, j, k} A^{(d^k n+1)f} T_{A^f, j,
k}^{\e_2} A^{(d^k n+2)f} \dots T_{A^f, j, k}^{\e_h} A^{(d^k
n+h)f}, \label{D:4}
\end{equation}
where $\e_1, \e_2 \dots, \e_h$, $h, n$ are defined as in
(\ref{D:2}).

It follows from definitions that the word $R_{A^f, j, k}$ is
conjugate (by $W^{-1}_{A^f, j, k}$) in rank $i-1$ to the word
 $w_{\Omega, k}( X_{A^f, j, k}, \bar Y_{A^f, j, k}) \overset {i-1}
\neq 1$.

The set $\Ss_i$ of defining words of rank $i$ consists of all
possible words $R_{A^f, j, k}$  given by
(\ref{D:3})--(\ref{D:4}) (over all equivalence classes
$\C_{A^f, j}(k)$, $A \in \X_i$, $k=1,2,\dots$). Finally, we put
$\R_i = \R_{i-1} \cup \Ss_i$ and set
\begin{gather}
G(i) = \la \A \; \| \; R=1, R \in  \R_i \ra  . \label{D:5}
\end{gather}
 The inductive definition of  groups $G(i)$, $i \ge 0$, is
 now complete and we can consider the limit group $G(\infty)$ given by
 defining words of all ranks $j =1,2, \dots$
\begin{gather}
G(\infty) = \la \A \; \| \; R=1, R \in \cup_{j=0}^\infty \R_j
\ra . \label{D:6}
\end{gather}

We will prove (in Lemma \ref{L5}) that $G(\infty)$ is the free
group of the variety $\MM_{\Omega}$ in the alphabet $\A$, that
is, $G(\infty)$ is naturally isomorphic to the quotient $F(\A)
/ W_{\Omega}(F(\A))$, where $W_{\Omega}(F(\A))$ is the verbal
subgroup of $F(\A)$ defined by the set $W_{\Omega} = \{
w_{\Omega, k}(x,y) \mid k=1,2, \dots \}$, and then show that
$G(\infty)$ is not abelian. But first we need to study the
presentation (\ref{D:5}) of $G(i)$. As in Sects. 29--30
\cite{O89}, following Lemmas \ref{L1}--\ref{L4} are proved by
induction on $i \ge 0$ (whose base for $i =0$ is trivial).

\begin{lemma} \label{L1}
The presentation $(\ref{D:5})$ of $G(i)$ satisfies the condition
$R$ of \cite[Sect. 25]{O89}.
\end{lemma}

\begin{proof} This proof is quite similar to the proof of Lemma 29.4
\cite{O89}.
Inductive references to Lemmas 30.3, 30.4, 30.5 \cite{O89} (in rank $i-1$)
are replaced by references to Lemma \ref{L4}. Note that, by Lemma \ref{L4}
and LPP, we have that
$$
|f|(d^kn+h)\le 100\zeta^{-1}(d^k n+h) <  d^{k+1} n
$$
(LPP: $\delta = h^{-1} \succ \zeta \succ \eta = d^{-1} \succ
\iota = n^{-1}$) which implies that, repeating the arguments of
Lemma 29.3  \cite{O89}, we can conclude that  the defining
relations $R$, $R'$ correspond to the same value of $k$ which
enables us to finish the proof of the analogue of Lemma 29.3 as
in  \cite{O89}.
\end{proof}

Now suppose that $X, Y$ are some words with $[X, Y] \overset
{i} \neq 1$ and $k$ is an arbitrary positive integer.
Conjugating the pair $(X, Y)$ in rank $i$ if necessary, we can
assume that $X \equiv B^{f_B}$,  $Y \equiv Z C^{f_C} Z^{-1}$,
where each of  $B, C$ is either simple in rank $i$ or a period
of some rank $\le i$ and, when $B^{f_B}$,  $C^{f_C}$ are fixed,
the word $Z$ is picked to have minimal length. Furthermore,
consider the following equalities
\begin{gather*}
[X^d, Y^d ] \overset {i} = W_D D^{f_D} W_D^{-1} ,  \\
[ [X^d, Y^d ]^d, X^{d q_k } ] \overset {i} = W_E E^{f_E}
W_E^{-1} ,
\end{gather*}
where each of $ D, E$ is either simple in rank $i$ or a period
of some rank $\le i$ and the conjugating words $W_D, W_E$ are
picked (when $D, E$ are fixed) to have minimal length. Without
loss of generality, we can also suppose that if $A_1 \in \{ B,
C, D, E\}$ is conjugate in rank $i$ to $A_2^{\pm 1}$, where
$A_2 \in \{ B, C, D, E\}$, then $A_1 \equiv A_2$.

\begin{lemma}  \label{L2}
In the foregoing notation, the following inequalities hold
\begin{gather}
0 < | f_D | \le 100 \zeta^{-1} ,  \label{L2:1} \\
\max (| B^{d f_B} |,  | C^{d f_C} | ) \le   \zeta^{-1} |
D^{f_D} | , \label{L2:2}
\\
|Z| < 15 \zeta^{-2} | D^{f_D} | , \label{L2:3} \\
|W_D| < 31 \zeta^{-2} | D^{f_D} | . \label{L2:4}
\end{gather}
\end{lemma}

\begin{proof} If $f_D = 0$, that is, $[X^d, Y^d ] \overset {i} = 1$ then,
by Lemmas 1,  25.2 and 25.12 \cite{O89}, we have $[X, Y ]
\overset {i} = 1$, contrary to the choice of $X, Y$. Hence $f_D
\neq 0$.

In view of equality $[ B^{d f_B} , Z C^{d f_C} Z^{-1}] \overset
{i} = W_D D^{f_D} W_D^{-1}$, there is a reduced diagram $\D$ of
rank $i$ on a thrice punctured sphere the labels of 3 cyclic
sections of whose boundary $\p \D$ are $B^{d f_B}$,  $B^{-d
f_B}$, $D^{f_D}$. If $| f_D | > 100 \zeta^{-1}$  then $\D$ is a
$G$-map (see Sect. 24.2 \cite{O89}) and, as in the proof of
Lemma 25.19 \cite{O89}, it follows from Lemma 24.8 \cite{O89}
that $ D^{f_D}  \overset {i} = 1$, contrary to $f_D \neq 0$.
Hence, $| f_D | \le 100 \zeta^{-1}$ and inequalities
(\ref{L2:1}) are proven.

If, say, $| D^{f_D} | < \zeta | B^{d f_B} |$, then
$\D$ is an $E$-map (see Sect. 24.2 \cite{O89}) and a contradiction to
$f_D \neq 0$  follows from Lemma 24.6 \cite{O89} exactly as above.
Hence, $| B^{d f_B} | \le  \zeta^{-1} | D^{f_D} |$
and inequalities (\ref{L2:2}) are proven.

It follows from definitions and Lemma 30.2 \cite{O89} that
$$
|Z| < 7 \zeta^{-1} (   | B^{d f_B} |+  | C^{d f_C} | + | D^{f_D} |  ) .
$$
Then, in view of inequalities (\ref{L2:2}), we have
$$
|Z| < 7 \zeta^{-1} ( 2\zeta^{-1}    +1  ) | D^{f_D} | <
15 \zeta^{-2} | D^{f_D} |  ,
$$
as claimed in (\ref{L2:3}).

By estimates (\ref{L2:2})--(\ref{L2:3}), we have
$$
|[ X^d, Y^d ]| = 2(  | B^{d f_B} |+  | C^{d f_C} | + 2 | Z |  )  <
2(2 \zeta^{-1}  + 30 \zeta^{-2} ) | D^{f_D} | =  61 \zeta^{-2} | D^{f_D} | .
$$
Hence, it follows from Lemmas \ref{L1} and 22.1 \cite{O89} that
$$
|W_D| <  (\gamma + \tfrac 1 2)( |[ X^d, Y^d ]| + | D^{f_D} | )
< 31 \zeta^{-2} | D^{f_D} |
$$
and Lemma \ref{L2} is proved.
\end{proof}

\begin{lemma} \label{L3}
In the foregoing notation, the following inequalities hold
\begin{gather}
0 < | f_E | \le 100 \zeta^{-1} , \label{L3:1} \\
| D^{d f_D} |  \le  \zeta^{-1} | E^{f_E} | , \label{L3:2} \\
| B^{d q_k  f_B } |  \le  \zeta^{-1} | E^{f_E} | , \label{L3:3} \\
|W_E| <  3  \zeta^{-1}  | E^{f_E} | . \label{L3:4}
\end{gather}
\end{lemma}

\begin{proof}  Assume that  $[ [ X^d, Y^d ]^d ,  X^{d q_k}] \overset {i} =
1$. Then, by Lemmas 1, 25.2  and 25.12 \cite{O89}, we have $[ [
X^d, Y^d ],  X ] \overset {i} = 1$ and so $[ Y^d X^d Y^{-d} ,
X^d ] \overset {i} = 1$. In view of Lemma 25.14 \cite{O89}, we
further  have $ [ X^d, Y^d ] \overset {i} = 1$.  Then, as
before, by Lemmas 1, 25.2, 25.12 \cite{O89}, we obtain that $ [
X, Y ] \overset {i} = 1$, contrary to the choice of $X$ and
$Y$.

It follows from definitions that
\begin{equation} [ W_D  D^{d
f_D } W_D^{-1},  B^{d q_k  f_B } ]  \overset {i} = W_E E^{f_E}
W_E^{-1} \label{L3:5}
\end{equation}
and so there is a reduced diagram  of rank $i$ on a thrice
punctured  sphere  the labels of 3 cyclic sections of whose
boundary are $D^{df_D}$, $D^{-df_D}$, $E^{f_E}$. Now we can
repeat proofs of inequalities (\ref{L2:1})--(\ref{L2:2}) to
obtain (\ref{L3:1})--(\ref{L3:2}).

In view of equality (\ref{L3:5}), there is a reduced diagram of
rank $i$ on a thrice punctured  sphere  the labels of 3 cyclic
sections of whose boundary are $B^{dq_k f_B }$,  $B^{-d q_k f_B
}$, $E^{f_E}$. Now we can see that the proof of inequality
\eqref{L3:3} is analogous to that of inequality \eqref{L3:2}.

As in the proof of Lemma \ref{L2}, we have from Lemmas \ref{L1}
and 22.1 \cite{O89} that
$$
| W_E | <   ( \gamma + \tfrac 12 )\cdot 2( 2| W_D | +
 | D^{d f_D} | + | B^{d q_k f_B} |  +  \tfrac 12  |  E^{f_E}  |) .
$$
Hence, by Lemma \ref{L2} and estimates (\ref{L3:2})--(\ref{L3:3}), we get
$$
| W_E | <   (1 + 2 \gamma ) (  (62 \zeta^{-2} d^{-1} +1 )| D^{d
f_D} | + | B^{d q_k f_B } | +   \tfrac 12  |  E^{f_E}  |) <   3
\zeta^{-1}| E^{f_E} |
$$
(LPP: $\gamma \succ \zeta \succ \eta = d^{-1}$)  and Lemma \ref{L3} is
proved.
\end{proof}

\begin{lemma} \label{L4}
Let $R_{A^f, j, k}$ be a defining word of rank $i+1$ defined by
$(\ref{D:3})$ if $\Omega(k) = 0$ or by   $(\ref{D:4})$ if
$\Omega(k) = 1$. Then $0 < | f| \le 100 \zeta^{-1}$, $|A| > d$,
the words $T_{A^f, j, k}$, $U_{A^f, j, k}$  do not belong to
the cyclic subgroup $\la A \ra$ of $G(i)$ and
$$
\max (|T_{A^f, j, k} | , |U_{A^f, j, k} |)  < d |A| .
$$
\end{lemma}

\begin{proof}  It follows from definitions that,
in the foregoing notation,  we can assume that
\begin{gather*}
A \equiv E ,  \\
T_{A^f, j, k}   \overset {i} = W_E^{-1} [B^{f_B}, Z C^{f_C} Z^{-1}] W_E ,
\\
U_{A^f, j, k}   \overset {i} = W_E^{-1} [B^{f_B}, Z C^{f_C} Z^{-1}]^2 W_E ,
\end{gather*}
and $f = f(A^f, j, k)$ is $f_E$. Hence, in view of Lemmas
\ref{L2} and \ref{L3}, we have that
\begin{gather*}
0 < | f| \le  100  \zeta^{-1} , \\
f_E^{-1}| E^{f_E} | =  |A| \ge 10^{-2} \zeta^2 | D^{d f_D} |
\ge 10^{-2} \zeta^3 d | B^{d f_B} | \ge 10^{-2} \zeta^3 d^2 > d
\end{gather*}
(LPP: $\zeta \succ \eta = d^{-1}$) and
\begin{multline*}
\max (|T_{A^f, j, k} | , |U_{A^f, j, k} |) \le 2 |W_E| + 8 |Z|
+ 4 | B^{f_B}| + 4 |
C^{f_C}| < \\
< ( 6 \zeta^{-1} + 120 \zeta^{-3} d^{-1} + 8 \zeta^{-2} d^{-2})
| E^{f_E} | < 7 \zeta^{-1} | E^{f_E} | < 700 \zeta^{-2} | E | <
d |A|
\end{multline*}
(LPP: $\zeta \succ \eta = d^{-1}$).

Assume that one of $T_{A^f, j, k}$, $U_{A^f, j, k}$  belongs to
$\la E \ra \subseteq G(i)$. Then one of $[B^{f_B}, Z C^{f_C}
Z^{-1}]$, $[B^{f_B}, Z C^{f_C} Z^{-1}]^2$ is  conjugate in rank
$i$ to  a power of $E$. However, by Lemmas \ref{L2} and
\ref{L3},
\begin{multline*}
\max ( |[B^{f_B}, Z C^{f_C} Z^{-1}]|, |[B^{f_B}, Z C^{f_C}
Z^{-1}]^2| ) <  (120
\zeta^{-3} d^{-1} + 8 \zeta^{-2} d^{-2}) | E^{f_E} | < \\
< 100 \zeta^{-1} (120 \zeta^{-3} d^{-1} + 8 \zeta^{-2} d^{-2})  | E | < |E|
\end{multline*}
(LPP: $\zeta \succ \eta = d^{-1}$),
whence $|T_{A^f, j, k}|, |U_{A^f, j, k}| < |E|$ which contradicts Lemmas
\ref{L1} and 25.17 \cite{O89}. Lemma \ref{L4} is proved.
\end{proof}

\begin{lemma} \label{L5}
The group $G(\infty)$, defined  by presentation $(\ref{D:6})$,
is naturally isomorphic to the free group $F(\A) /
W_{\Omega}(F(\A))$ of the variety $\MM_{\Omega}$ in the
alphabet $\A$.
\end{lemma}

\begin{proof} It follows from the  definition  of defining words
of the group $G(\infty)$ that each of them is in
$W_{\Omega}(F(\A))$  and so there is a natural epimorphism
$$
G(\infty) \to F(\A) / W_{\Omega}(F(\A)) .
$$

Suppose that $\wtl X$, $\wtl Y$ are some words in $\A^{\pm 1}$ and
\begin{equation}
\label{L5:1} w_{\Omega, k}(\wtl X, \wtl Y ) \neq  1
\end{equation}
in $G(\infty)$ for some  integer $k >0$. Let $A$ be a period of
some rank such that $A^f$ for some $f$ is conjugate in
$G(\infty)$ to $v_{k}(\wtl X, \wtl Y )$. (The existence of such
an $A$ follows from definitions; see also Lemma 18.1
\cite{O89}.)  Note that, in view of (\ref{L5:1}), $[\wtl X,
\wtl Y ]   \neq 1$ in $G(\infty)$. Hence, by Lemmas \ref{L2}
and \ref{L3}, we can replace the pair $(\wtl X, \wtl Y )$ by a
conjugate in the group $G(\infty)$ pair $(X, Y)$ such that $X
\equiv B^{k_B}$, $Y \equiv Z C^{k_C} Z^{-1}$, and $v_{k}(X, Y )
= W_A A^f W_A^{-1}$ in $G(\infty)$, where $B$, $C$ are some
periods, $|f|>0$,
\begin{multline*}
|X^d | + |Y^d | = | B^{d k_B} | + |C^{d k_C} | +2 |Z| <
(2\zeta^{-1} + 30 \zeta^{-2}) |D^{k_D} | <  31 \zeta^{-3}
d^{-1} |A^f| ,
\end{multline*}
and $|X^{d q_k}|\le \zeta^{-1}|A^f|$.

Hence,
\begin{multline*}
| v_{k}( X, Y ) |  \le 2(2d( |X^d| + |Y^d|) + |X^{d q_k}|) < \\
< 2(62 \zeta^{-3} + \zeta^{-1}) |A^f|< 10^3 \zeta^{-3} |A^f| \le 10^5
\zeta^{-4} |A|
\end{multline*}
for  $0 < |f| \le 100 \zeta^{-1} $ by Lemma \ref{L3}. Consider
a reduced annular diagram $\D$ of some rank $i'$ for conjugacy
of $v_{k}( X, Y )$  and $A^f$. By Lemmas \ref{L1} and 22.1
\cite{O89}, $\D$ can be cut into a simply connected  diagram
$\D_1$ along a simple path $t$ which connects points on
distinct components of $\p \D$ with $| t| < \gamma |\p \D |$.
Therefore,
$$
|\p \D_1 | < (1 + 2 \gamma) |\p \D | < (1 + 2 \gamma)(10^5 \zeta^{-4}+
100\zeta^{-1} )|A|
< \tfrac 12 n |A|
$$
(LPP: $\gamma \succ \zeta \succ \eta = d^{-1} \succ \iota =
n^{-1}$). Then, by Lemmas \ref{L1}, 20.4 and 23.16 \cite{O89}
applied to $\D_1$, the diagram $\D_1$ contains no 2-cells of
rank $ \ge |A|$, whence $\D_1$, $\D$ are diagrams of rank $|A|
-1$. Since $A \in \X_{|A|}$, it follows from the construction
of defining words of rank $|A|$ that there will be a defining
word in $\Ss_{|A|}$ which guarantees that $w_{\Omega, k}(X, Y)
\overset {|A|} = 1$.  A contradiction to assumption
(\ref{L5:1}) proves that $G(\infty)$ is in $\MM_{\Omega}$ and
Lemma \ref{L5} is proved.
\end{proof}

{\em Proof of Theorem.} By Lemma \ref{L5}, $G(\infty)$ is the
free group of $\MM_{\Omega}$ in $\A$. Assume that $[a_1,
a_2]=1$ in $G(\infty)$. Then there exists an $i>0$ such that
$[a_1, a_2] \overset {i} =1$. This, however, contradicts Lemmas
\ref{L1} and 23.16 \cite{O89}. Thus, $G(\infty)$ is not abelian
and Theorem is proved. \qed

\end{document}